\documentclass[a4paper,11 pt]{amsart}
\usepackage{srollens-en}[2008/11/11]
\usepackage{epsfig}

\newcounter{saveenum}

\title{Compact Moduli for certain Kodaira fibrations}

\author{S\"onke Rollenske}

\address{Dr. S\"onke Rollenske\\
Mathematisches Institut \\
 Rheinische Fried\-rich{\-}-Wil\-helms-Uni\-versi\-t\"at Bonn \\
 Endenicher Allee 60  \\
 53115 Bonn, 
Germany}
\email{srollens@math.uni-bonn.de}

\DeclareFontFamily{OT1}{rsfs}{}
\DeclareFontShape{OT1}{rsfs}{n}{it}{<-> rsfs10}{}
\DeclareMathAlphabet{\curly}{OT1}{rsfs}{n}{it}

\DeclareMathOperator{\Def}{{Def}}
\newcommand{\T}{{\mathbf T}}
\newcommand{\beq}[1]{\begin{equation}\label{#1}}
\newcommand{\eeq}{\end{equation}}

\newcommand{\CC}{{C\times C}}
\newcommand{\MKF}{{\gothM}^{KF}}
\newcommand{\MKFbar}{\overline{\gothM}^{KF}}
\newcommand{\M}[1]{{\overline{\gothM}^{#1}}}

{\theoremstyle{Lehn-it}
\newtheorem{ass}[theo]{Assumptions}}

\begin{document}

\begin{abstract}
We describe explicitly the possible degenerations of a class of double Kodaira fibrations in the moduli space of stable surfaces.
Using deformation theory we also show that under some assumptions we get a connected component of the moduli space of stable surfaces.
\end{abstract}

\subjclass[2000]{14J29; 14J10, 14D20}

\maketitle

\section*{Introduction}
It is a general fact that moduli spaces of \emph{nice} objects in algebraic geometry, say smooth varieties,  are often non-compact. But usually there is a modular compactification where the boundary points correspond to related but more complicated objects.

Such a modular compactification has been known for the moduli space $\km_g$ of smooth curves of genus $g$ for a long time and in \cite{ksb88} Koll\'ar and Shepherd-Barron made the first step towards the construction of a modular compactification $\M{}$ for the moduli space $\mathfrak M$ of surfaces of general type via so called stable surfaces; the boundary points arise from a stable reduction procedure. 

But even 20 years later very few explicit descriptions of compact components of $\M{}$ have been published. The main idea in all approaches is to relate the component of the moduli space one wishes to study to some other moduli space, where a suitable compactification is known. Products of curves and surfaces isogenous to a product of curves have been treated by van Opstall \cite{vop05, vop06} and a recent paper of Alexeev and Pardini \cite{ale-pard09} studies Burniat and Campedelli surfaces relating them to hyperplane arrangements in (a blow-up of) $\IP^2$.

The aim of this paper is to study the irreducible resp. connected components of the moduli space of stable surfaces containing very simple Galois double Kodaira fibrations (see Section \ref{kodfib} for the precise definition) and we do this in 2 steps: first we give an explicit description of the stable degenerations  and then we study their deformations to show if we get  connected components of $\M{}$.

The starting point of our study is a joint paper with Fabrizio Catanese where we showed that the moduli space of such (and more general) surfaces can be identified with the moduli space of certain curves with automorphisms and yields connected components of the moduli space of surfaces of general type.

The surfaces we are interested in  are ramified covers $\psi:S\to C_1\times C_2$  of products of curves and we have precise control over the branch divisor $B$. This enables us to perform the stable reduction procedure explicitly, first for the pair $(C_1\times C_2, B)$ in Section \ref{pairs} and then for the Kodaira fibrations themselves in Section \ref{fibs}. It turns out that  degenerations that occur are local complete intersections and their normalisation is smooth (Theorem \ref{main}).

In Section \ref{def} we use deformation theory to study the scheme structure of the moduli space.   We show that  the deformations of standard Kodaira fibrations are unobstructed and are exactly the ones described in \cite{cat-roll06} if some cohomology groups naturally associated to the covering $\psi$ vanish (Theorem \ref{defo1}). Under similar assumptions we are able to control all deformations of the degenerations (Theorem \ref{main2}).
In a special case the assumptions are easy to check and we get
\begin{custom}[Corollary \ref{cyclic}]
If $\psi:X\to C\times C$ is a smooth very simple Kodaira fibration such that $\psi$ is  a cyclic covering and let $\gothN$ be the irreducible component of the moduli space of surfaces of general type containing (the class of) $X$. 
Then the closure of the component $\overline\gothN \subset \M{}$ is a connected component of $\M{}$.
\end{custom}

A rather different and less explicit approach to the construction of a compact moduli space for fibred surfaces has been described by Abramovich and Vistoli \cite{av2000}.

\subsection*{Acknowledgements:} Numerous people answered questions concerning particular parts of this article, including Paolo Cascini, Alessio Corti, Donatella Iacono, S\'andor Kov\'acs, Michael L\"onne and Thomas Peternell. Special thanks goes to Fabrizio Catanese for suggesting Lemma \ref{parallel}. A request of the referee encouraged me to push through the calculations in Section \ref{new} that lead to Theorem \ref{main2}. Part of this work was written during a stay at Imperial College London supported by a Forschungsstipendium of the DFG. The author was also supported by the Hausdorff Centre for Mathematics in Bonn. The DFG-Forschergruppe ``Classification of algebraic surfaces and compact complex manifolds'' made a visit to Bayreuth possible.

\section{Preparations}
\subsection{Stable surfaces and some other moduli spaces}
We will start this section by listing some (coarse) moduli spaces that we will use in the sequel, mainly to fix the notation, and give references to where a construction and more information can be found. The moduli space of (smoothable) stable surfaces will be discussed a bit more in detail.
\begin{itemize}
 \item Let $\km_g$ be the moduli space of smooth projective curves of genus $g$ and $\overline{\km}_g$ the moduli space stable curve (see e.g. \cite{Harris-Morrison}). 
\item The corresponding moduli spaces $\km_g(G)$ and $\overline{\km}_g(G)$ parametrising smooth resp. stable curves together with a fixed group of automorphisms. The moduli space $\overline{\km}_g(G)$ is finite over a closed subvariety of the moduli space of curves and the tangent space at a point $[C]$ is $\Ext^1(\Omega_C, \ko_C)^G$ (see \cite{tuff93} and also \cite[Prop. 2.9]{vop06}).
\item The quasi-projective coarse moduli space $\mathfrak M_{a,b}$ of canonically polarised surfaces of general type with canonical singularities and fixed invariants $a=K_S^2, b=\chi(\ko_S)$ constructed by Gieseker \cite{gieseker77}. We denote by $\mathfrak M$ the disjoint union of all $\mathfrak M_{a,b}$.
\item The moduli space of  stable surfaces $\M{}=\M{st}$, which contains $\mathfrak M$ as an open (but not dense) subset. 
Its construction goes back to Koll\'ar and Shepherd-Barron \cite{ksb88}, and we will give some more details below. 
We denote by $\M{sm}$ the closure of $\gothM$ in $\M{}$ and call it the moduli space of smoothable surfaces of general type.
\end{itemize}
Note that, strictly speaking, a surface in $\M{sm}$ need not to be smoothable in the ordinary sense -- we only ask that some small deformation has canonical singularities.

In analogy to the case of curves we need to define a class of singular surfaces that is big enough to allow for compact moduli spaces. It is convenient to first recall the definition of log-canonical singularity. We will need this notion also for pairs.
\begin{defin}
Let $X$ be a normal surface and $B\subset X$ a (possibly empty) $\IQ$-Weil-divisor such that the log-canonical divisor $K_X+B$ is $\IQ$-Cartier. Let $\pi:\tilde X\to X$ be a log-resolution of singularities. In other words, $\tilde X$ is smooth and denoting the strict transform of $B$ with $\tilde B$ and the exceptional divisor with $E=\sum_i E_i$, the sum $\tilde B +E$ is a global normal crossing divisor. Then the pair $(X, B)$ is called log-canonical (lc) resp. canonical if in the expression
\[ K_{\tilde X}+\tilde B \equiv \pi^*(K_X+B) + \sum a_i E_i\]
all $a_i\geq -1$ resp. $a_i\geq 0$.
\end{defin}
Canonical surface singularities without boundary are exactly rational double points. The notion we are aiming at is some kind of non-normal analog of log-canonical singularities.
\begin{defin}
 Let $S$ be a projective surface and $B=\sum b_i B_i$ an effective $\IQ$-Weil divisor on $S$ with coefficients $0<b_j\leq 1$.

The pair $(S,B)$ is said to have slc singularities if
\begin{enumerate}
 \item $S$ is Cohen-Macaulay,
\item $S$ has at most normal crossing singularities in codimension 1 and $B$ does not contain any component of the normal crossing locus,
\item $K_S+B$ is $\IQ$-Cartier,
\item denoting by $\nu:X^\nu\to X$ the normalisation, the pair 
\[(X^\nu, (\text{double locus})+\inverse\nu B)\]
 is log canonical. By double locus we mean the preimage of the 1-dimensional part of $X_{\text{sing}}$, which outside a finite number of points coincides with the normal crossing locus.
\end{enumerate}
The pair $(S,B)$ is called stable if it has slc singularities and the $\IQ$-line bundle $K_S+B$ is ample.
\end{defin}
The original definition was posed in \cite[Section 4]{ksb88} where one can also find a classification in the case $B=0$.

In the construction of the moduli space, especially if one wants to parametrise pairs,  several technical issues arise, see \cite{kovacs05} or \cite{alexeev08} for an overview.  In particular, it turns out that one has to restrict the families that are allowed in the moduli-functor if one wants the basic invariants to remain constant in a family. As usual, we define the reflexive powers of a sheaf $\kf$ by $\kf^{[n]}:=(\kf^{\tensor n})^{**}$.
\begin{defin}\label{allowablefamily}
A morphism $f:(\kx,B)\to \Delta$ is called a \defobj{weakly stable  family of stable surfaces} if
\begin{enumerate}
 \item $f$ and $f\restr B$ are flat and projective and $f$ has connected fibres,
\item $\omega_{\kx/\Delta}+B$ is a relatively ample $\IQ$-line bundle,
\item For all $t\in \Delta$, the fibre $\kx_t$ is a stable surface.
\setcounter{saveenum}{\value{enumi}}
\end{enumerate}
If $B=0$ and in addition we have
\begin{enumerate}
\setcounter{enumi}{\value{saveenum}}
\item (Koll\'ar's condition) For all $t\in \Delta$, $k\in \IZ$, taking reflexive powers commutes with base-change, that is $\omega_{\kx/\Delta}^{[k]}\restr{\kx_t}\isom \omega_{\kx_t}^{[k]}$.
\end{enumerate}
then $f:\kx\to \Delta$ is called \defobj{admissible family of stable surfaces}.
\end{defin}
Since we are interested only in  degenerations of smooth surfaces, we will usually assume that $\Delta$ is a smooth curve and that the general fibre of $f$ is canonical.

A recent construction of $\M{}$ using stacks can be found in \cite{abr-hass09}, including a proof that $\M{sm}$ is projective. The bigger moduli space of stable surfaces should also be projective but I do not know a suitable reference for this. In particular  stable reduction  as explained in the next section seems not to be known for arbitrary admissible families of stable surfaces.

\subsection{Construction of the boundary points}\label{stabred}

We need to describe how to obtain surfaces corresponding to the boundary points in $\M{sm}\setminus\gothM$ as degenerations of smooth surfaces.  Let $\Delta$ be the unit disc and $\Delta^*=\Delta\setminus\{0\}$. We will need the construction also for of families of log-surfaces.

Suppose we have a family of log-surfaces $f^0:(\kx^0,B^0)\to \Delta^*$, that is, both $f^0$ and $f^0\restr{B^0}$ are flat, projective maps and for each $t\in \Delta^*$ the fibre $(\kx_t^0, B^0_t)$ is a log-surface. Suppose in addition that
\begin{itemize}
 \item all fibres $\kx_t^0$ are canonical,
\item $K_{\kx^0/\Delta^*}+B^0$ is a relatively ample $\IQ$-line bundle.
\end{itemize}
 Then one constructs the degeneration of this family in the following steps (see  \cite[Theorem 7.62]{Kollar-Mori}):
\begin{enumerate}
 \item Choose any extension of $f^0$ to a projective morphism $f:(\bar \kx,\bar B)\to \Delta$.
\item Apply the semi-stable reduction theorem to a log-resolution of $(\bar \kx,\bar B)$ obtaining, possibly after a finite pullback ramified only over the central fibre and normalisation, a family of log-surfaces $(\tilde \kx, \tilde B)\to \Delta$ such that the total space is smooth, $(\tilde \kx, \tilde B+\tilde \kx_0)$ is log canonical and $K_{\tilde \kx/\Delta}+\tilde B$ is relatively big over $\Delta$.
\item Now let $f:(\kx,B)\to \Delta$ be the relative log-canonical model of $(\tilde \kx, \tilde B)$, whose existence is guaranteed by the log-minimal model program.
\end{enumerate}
The resulting family is, due to the possible pullback, not unique but the central fibre is; we call it the stable degeneration of the family $f^0$.

If some of the coefficients of $B$ are strictly smaller than 1 then there are examples due to Hassett \cite[Example 5.1]{alexeev08} which show that $B_0$ can have embedded points so we do not obtain a family of log-surfaces. These problems do not occur, if $B$ is an integral divisor  \cite{hassett01} or if there is no boundary. More general results have been obtained by Alexeev \cite{alexeev08}.
%

\subsection{Very simple Galois Kodaira Fibrations}\label{kodfib}
We will now introduce the class of surfaces of general type that we are interested in. In this section all surfaces and curves are smooth.

\begin{defin}\label{defSKF} A \emph{Kodaira fibration} is a smooth
fibration $\psi_1:S\to C_1$ of a compact complex surface over a
compact complex curve, which is not a holomorphic fibre bundle.

$S$ is called a \emph{double Kodaira fibred surface} if it admits
a \emph{double Kodaira fibration}, i.e.,  a surjective holomorphic map
$ \psi : S \to C_1 \times C_2$ yielding two  Kodaira fibrations
$\psi_i : S \to C_i$ $(i=1,2)$.

    Let $B\subset C_1\times C_2$ be the branch divisor of $\psi$.  If
$B$ is smooth and both  projections $pr_{C_j}\restr{B}:B\to C_j$ are 
\'etale we call
$ \psi : S \to C_1 \times C_2$ a
\emph{double \'etale Kodaira fibration}.

If the ramified cover $\psi$ is Galois, i.e., the quotient map for the action of a finite group, we call $S$ a Galois double Kodaira fibration.
\end{defin}

It is not difficult to see that if $S\to C_1\times C_2$ is  a double Kodaira surface then the genus of $C_i$ is at least 2 and thus $S$ is a surface of general type.
A situation in which the branch divisor is particularly easy to handle is the following.
\begin{defin}
A double \'etale Kodaira fibration $\psi: S\to C_1\times C_2$ is called very simple if $C_1=C_2=C$ and the branch divisor is a disjoint union of graphs of automorphisms of  $C$. It is called standard, if there is an \'etale map $\phi:C\times C \to C_1\times C_2$ such that $\phi^*S\to C\times C$ is very simple.
\end{defin}
\begin{rem}
It is yet unclear whether every double \'etale Kodaira fibration is standard: let $B\subset C_1\times C_2$  be a divisor in a product of 2 curves mapping \'etale to both sides. By taking $\tilde C_1\to C_1$ to be a Galois cover dominating all components of $B$ we obtain after pullback that $\tilde B\subset \tilde C_1\times C_2$ is composed of graphs of \'etale maps. But it is not at all clear that after further pullback we can arrive at graphs of automorphisms. On the other hand we do not know of an explicit example where this is not possible. 
\end{rem}
In \cite{cat-roll06} we described an effective method of construction for Galois double Kodaira fibrations. Essentially, given two curve $C_1, C_2$ and  a branch divisor $B\subset C_1\times C_2$ mapping \'etale to both curves we can construct plenty of Galois double Kodaira fibration (after finite \'etale pullback). This generalises a classical construction used by Kodaira and Atiyah to give examples of fibre bundles where the signature is not mutliplicative (see \cite{BHPV}, Section V.14). The basic idea is as follows: assume that we have a curve $C$ of genus at least 2 and a group $H$ acting on $C$ without fixed points. Then the graphs of the automorphisms do not intersect and $B=\bigcup_{\phi\in H} \Gamma_\phi$ is a smooth divisor in $C\times C$. After a suitable pullback $\pi:\tilde C\times C\to C\times C$ the divisor $\pi^*B$ will be divisible in $\mathrm{Pic}(\tilde C\times C)$ and we get a Kodaira fibration as a cyclic covering. With a further pullback we can make the branch divisor into a union of graphs of automorphisms again, obtaining a very simple Kodaira fibration.

We were able to give a very explicit description of the moduli space of such surfaces.
Let $\MKF\subset \gothM$ be the subset of the moduli space of surfaces of general type containing very simple Galois Kodaira fibrations and $\MKFbar$ its closure in the moduli space of stable surface. If $S$ is a very simple  Kodaira fibration then we denote by $\MKF_S$ the connected component of $\MKF$ containing (the class of) $S$.
\begin{theo}[\cite{cat-roll06}, Theorem 6.5]\label{moduliset}
Let $\psi: S\to C\times C$ be a very simple Galois Kodaira fibration, $\ks \subset \Aut(C)$ such that the branch divisor $B=\sum_{\sigma\in \ks} \Gamma_\sigma$. Let $H$ be the subgroup of $\Aut(C)$ generated by $\ks$. 
Then $\MKF_S$ is a connected component of the moduli space of surfaces of general type and  there is a natural map  $\km_{g(C)}(H)\to\MKF_S$ that  is an isomorphism on geometric points.
\end{theo}
In other words, any actual family of curves with the prescribed automorphism group gives rise to a family of Kodaira fibrations but we cannot detect obstructed deformations, or additional automorphisms of $S$ that do not preserve the map $\psi$.
The original theorem is formulated for standard Kodaira fibrations but we will only need this simple form. The issue of the scheme structure of $\MKF$ will be addressed in Section \ref{def}.

\section{Degenerations of very simple Galois Kodaira fibrations}
\subsection{Quotient families}
We start with some general considerations. Let $f:X\to Y$ be a Galois cover with ramification divisor $R\subset X$ and branch divisor $B\subset Y$. By the Hurwitz formula we have  $K_X=f^*K_Y+\sum_i (\nu_i-1) R_i$ where $\nu_i$ is the ramification order along $R_i$. Since the covering is Galois two components of $R$ that map to the same component of $B$ have the same ramification order and thus we can write 
$ K_X=f^*(K_Y+\sum_j\frac{\nu_j-1}{\nu_j}B_j)$. This allows us to compare the numerical properties of $K_X$ and the $\IQ$-divisor $K_Y+\sum_j\frac{\nu_j-1}{\nu_j}B_j$.
\begin{lem}\label{quotlem}
Let $\Delta$ be a smooth curve, $f:\kx\to \Delta$ be a flat, projective family of surfaces together with an action of a finite group $G$ preserving the fibres. Consider the quotient family 
\[ \xymatrix{ \kx \ar[rr]^\pi\ar[dr]_f && \kx\slash G =\ky \ar[dl]^g\\& \Delta}\]
and let $\kb\subset \ky$ be the branch divisor, that is, the divisorial part of the branch locus, with the appropriate multiplicities such that $K_{\kx\slash \Delta}=\pi^*(K_{\ky/\Delta}+\kb)$.

If $\kx\to \Delta$ is a weakly stable  family of stable surfaces then so is $(\ky,\kb)\to \Delta$.
\end{lem}
\begin{proof}
Since $G$ acts fibrewise, every irreducible component of $\kx$ dominates $\Delta$ if and only if every irreducible component of $\ky$ does and thus, by \cite[Proposition III.9.7]{Hartshorne}, $f$ is flat if and only if $g$ is. The same is true for projectivity.

Let $\pi_t:\kx_t \to \ky_t$ be the restriction of $\pi$ to some fibre, $x\in \kx_t$ and $y=\pi_t(x)$. By construction we get an inclusion of local rings $\ko_{\ky_t,y}\into \ko_{\kx_t,x}$. Indeed, if $G_x$ is the stabiliser of $x$ in $G$ then $\ko_{\ky_t,y}=\ko_{\kx_t,x}^{G_x}$. Since we assumed $\kx_t$ to be reduced and Cohen-Macaulay the same holds for $\ky_t$ where the second property is proved via an averaging argument for a regular sequence.

In order to prove that $g$ is a weakly stable family we also  need to show that the branch divisor $\kb$ is flat over $\Delta$ which amounts to the fact that $\pi$ is not ramified along any irreducible component of any fibre. But for every $t\in \Delta$ the map $\kx_t\to \ky_t$ is flat by base-change and since all fibres of $f$ are reduced it can not be ramified along any irreducible component by generic smoothness.

By \cite[Lemma 4.3]{ale-pard09} the fibre $\kx_t$ is slc if and only if the pair $(\ky_t, \kb_t)$ is slc and in particular it makes sense to ask if  $K_{\ky\slash \Delta}+\kb$ is relatively ample. This follows from the Nakai-Moishezon criterion: if $C$ is a curve contained in a fibre of $g$ then 
\[ (K_{\ky\slash \Delta}+\kb).C=\frac{1}{|G|}\pi^*(K_{\ky\slash \Delta}+\kb).\pi^*C=\frac{1}{|G|}K_{\kx/\Delta}.\pi^*C>0\]
because $K_{\kx/\Delta}$ is relatively ample. Also $(K_{\ky_t}+\kb_t)^2=1/{|G|}K_{\kx_t}^2>0$ and we see that $K_{\ky\slash \Delta}+\kb$ is relatively ample as well which concludes the proof.
\end{proof}
The above Lemma enables us to relate degenerations of very simple Kodaira fibrations to a situation we control better.
\begin{prop}\label{lcmod} Let $\Delta^*$ be the pointed disk and let $f^0:\kx^0\to \Delta^*$ be an admissible family of Galois double Kodaira fibrations with Galois group $G$. Then, possibly after a finite pullback, there are two families of curves $\kc_i\to \Delta^*$ fitting in the diagram
\[ \xymatrix{ \kx^0 \ar[rr]^\pi\ar[dr]_f && \kx^0\slash G =\kc_1\times \kc_2 \ar[dl]^g\\& \Delta^*}.\]
Denoting by $\kb^0\subset \kc_1\times \kc_2$ the branch divisor of $\pi$ (with the appropriate multiplicities) the stable degeneration of $f$ is a ramified cover of the stable degeneration of the family of log-surfaces  $g:(\kc_1\times \kc_2, \kb^0)\to \Delta^*$.
\end{prop}
\begin{proof} Since the Galois group $G$ is finite and its action on the differentiable manifold underlying the fibres $\kx_t$ $(t\neq 0)$ is fixed, possibly after a finite \'etale pullback, the action of $G$ on the fibres extends to a fibrewise action of $G$ on the whole family $\kx^0$ such that the quotient is a family of products of curves $g^0: \kc_1\times \kc_2 \to \Delta^*$.

Now let $f:\kx \to \Delta$ be a stable degeneration of $f^0$. Using the same argument as in the proof of the separatedness of the moduli space (see e.g. \cite[Corollary 5.15]{kovacs05}) we see that the action of $G$ on $\kx^0$ extends to the whole family $\kx$. By Lemma \ref{quotlem} the quotient together with the branch divisor  is a weakly stable family $g:(\kx/G, \kb)\to \Delta$ of stable log-surfaces. 
But by construction the families $g$ and $g^0$ coincide (up to finite \'etale pullback) over $\Delta^*$ and we get the claimed result because the stable degeneration is unique by separatedness. \end{proof}

\subsection{Stable reduction of the family of pairs}\label{pairs}
We have seen above that in order to understand slc degenerations of very simple Galois Kodaira fibrations we need to understand the stable reduction of families of the following type.
\begin{ass}\label{assum}
Let $f:\kc\to \Delta$ be a family of stable curves over the unit disk with smooth general fibre such that a group of 
automorphisms $H$ acts fibrewise on the fibration $f$. Let $\ks \subset H$ be a subset and let \[\kb=\sum_{\phi\in\ks}\lambda_\phi\Gamma_\phi\subset \ky:=\kc\times_\Delta \kc\]
be the union of the (fibrewise) graphs of the automorphisms in $\ks$ with some rational coefficients $0<\lambda_\phi<1$. We further assume that  on the general fibre the graphs of two different automorphisms do not intersect.
\end{ass}

We will now analyse the local structure of $(\ky,\kb)$ and see that while it is not stable itself (as soon as the central fibre is singular)
only a simple modification is needed to stabilise it. Since the general fibre $( \ky_t, \kb_t)$ consists of a smooth surface containing a smooth $\IQ$-divisor we only have to consider the central fibre $(\ky_0, \kb_0)$.

If the boundary is empty then the central fibre is slc, the only singularities except normal crossings being degenerate cusps:
%

\begin{lem}\label{prodnode}
Consider the degenerate cusp surface singularity $0\in Z=\spec\IC[x_1, \dots, x_4]/(x_1x_2, x_3x_4)$ given by the product of two 1-dimensional nodes. Let $A$ be a curve passing through $0$ but not contained in the double locus. Then the pair $(Z,\lambda A)$ is slc if and only if $\lambda=0$.
\end{lem}
\begin{proof} We work in the local model described above and consider  the blow-up  $\pi:\tilde Z\to Z$ in $0$ inside $\IC^4\times \IP^3$. Since $Z$ is a cone the exceptional divisor $E$ a cycle of 4 lines in $\IP^3$. Note that $\tilde Z$ is in fact semi-smooth and a good semi-resolution of $Z$ in the sense of \cite{ksb88}.

We can choose $A$ to be given by the ideal $(x_1, x_2, x_3-x_4)$ and see that its strict transform $\tilde A$ intersects $E$ in a single point $(0,(0:0:1:-1))$ and that $\pi^*A=\tilde A+E$. Let us now calculate discrepancies. 

By the adjunction formula we have
\[0=(K_{\tilde Z} +E).E=(\pi^* K_Z +a_0E+E).E=(a_0+1)E^2,\]
thus $a_0=-1$,  $K_{\tilde Z}=\pi^*K_Z -E$ and the singularity is slc if there is no boundary. 
If we add the boundary divisor we see that 
\[ K_{\tilde Z} + \lambda \tilde A = \pi^*(K_Z+\lambda A)-(1+\lambda )E \]
and the pair is not slc if $\lambda\neq 0$.\end{proof}

It turns out that these are the only kind of non-slc points that can occur in the chosen degeneration:
\begin{lem}\label{parallel}
In the situation of \ref{assum} assume that we have two automorphisms $\phi\neq\phi'\in \ks$ such that their graphs $A:=\Gamma_{\phi}$ and $A':=\Gamma_{\phi'}$ intersect in a point $x=(p,q)\in C_0\times \kc_0=\ky_0$. Then
 $p$ and $q$ are both nodes and $(\ky_0, \lambda A +\lambda ' A' )$ is slc at $x$ if and only if $\lambda=\lambda'=0$. 

The automorphism $\inverse{\phi'}\circ \phi$ preserves the local branches at $p$ and the divisors $A$ and $A'$ intersect transversely on each irreducible component of the normalisation of $\ky_0$.
\end{lem}
For the last statement it is actually enough to look at the irreducible components of a small analytic neighbourhood of $(p,q)$.
\begin{proof}
By composing both automorphisms with $\inverse{\phi'}$ we can assume that $\phi'=id$ and $p=q$. In other words, the automorphism $\phi$ has a fixed point at $p$. We first have to show that $p$ cannot be a smooth point of $\kc_0$. Indeed, in that case the threefold $\ky$ is smooth in a neighbourhood of $(p,p)$ and the two divisors $\Gamma_\phi$ and $\Gamma_{id}$ intersect in a curve which is completely contained in the central fibre $\ky_0$ by our assumptions. In other word $\phi$ acts as the identity on some irreducible component $\kc_0'$ of $\kc_0$. Let $p$ be a node where $\kc_0'$ meets another component. Locally around $p$ the total space is described by the equation $xy-t^m=0$ in $\IC^3$ where $t$ is the coordinate on the base $\Delta$ of the fibration. The (linearised) action of $\phi$ is trivial on one component, say $t=y=0$, and preserves the fibres and therefore $\phi^*$ acts trivially on both $x$ and $t$. Since it also has to preserve the equation it necessarily acts trivially on $y$ and thus also on the other component. Since the curve is connected $\phi=id_{\kc_0}$. But this is a contradiction since a non-trivial automorphism cannot degenerate to the identity and thus the graphs cannot meet in a smooth point of $\ky_0$.

So assume now that $p$ is a node of $\kc_0$ and that $\phi$ interchanges the local branches at $p$. In the local model as above the (linearised) action is given by 
\[ (x,y,t)\mapsto (\zeta y, \inverse\zeta x, t)\]
for some root of unity $\zeta\neq 1$ because the group $H$ acts fibrewise and thus leaves $t$ fixed while preserving the equation. But for small $t\neq 0$ there is a non-trivial solution $x$ to $\inverse\zeta x^2=t^m$ and the point $(x, \inverse\zeta x,\inverse\zeta x^2)$ is a fixed point for the action of $\phi$ which contradicts the assumption that the divisors are disjoint in the general fibre. Thus $\phi=\inverse{\phi'}\circ \phi$ preserves the local branches at $p$.

 When restricted to one of the irreducible components (either in a neighbourhood of $(p,p)$ or on the normalisation) we are locally looking at the intersection of divisors given by $x-y$ and $x-\zeta y$ in $\IC^2$ and thus their intersection is transverse. The statement on the singularity not being slc has been proved in Lemma \ref{prodnode}.\end{proof}

We will now describe the stable reduction procedure (see Section \ref{stabred}) in the local model  where $\ky$ is given by $(x_1x_2-t, x_3x_4-t)\subset \IC[x_1,\dots,x_4,t]$; the total space of $\kc$ is smooth and the threefold $\ky$ has an isolated singularity at the point $x=(0,0,0,0,0)$. Blowing up this singular point we obtain a smooth threefold $\tilde \ky\to \ky$. The central fibre is a normal crossing divisor, union of the blow-up of $\ky_0$ described in the proof of Lemma \ref{prodnode} and the exceptional divisor $Q\isom \IP^1\times \IP^1$, which are glued together along a the exceptional cycle of four lines as depicted in Figure 1.

\begin{figure}[h]\caption{Stable reduction at a degenerate cusp where 2 components of $\kb_0$ meet.}\label{stabredfig}
\centering{\epsfig{file=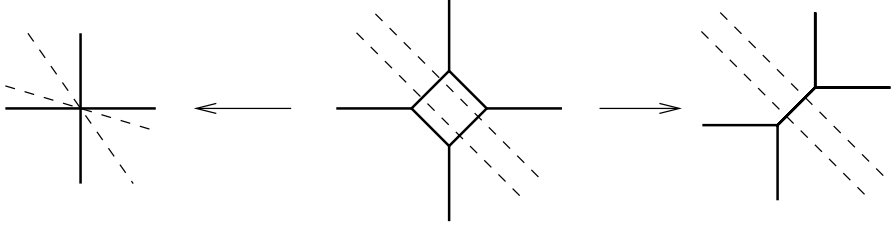}} 
\end{figure}

We may assume that all components of  $\kb_0=A_1+\dots +A_k$  pass through $x$. By Lemma \ref{parallel} all corrresponding automorphisms preserve the local branches at $x$ and thus $\kb_0$ is contained in the union of the components $t=x_1=x_3=0$ and $t=x_2=x_4=0$. When we only look at one irreducible component the $A_i$ meet transversely at $x$; all branches get separated after the blow-up and the strict transform $\tilde A$ intersect the double locus transversely; on the central component $Q$ we get several parallel lines.

Denoting by $\tilde \ky_0$ the central fibre of $\tilde \ky$  we see that the pair $(\tilde \ky_0, \tilde \kb)$ is slc and thus by inversion of adjunction $(\tilde\ky, \tilde \kb+\tilde \ky_0)$ is lc as required.  In order to to take the relative canonical model over $\Delta$ we have to see on which curves in $\tilde \ky_0$ the log-canonical divisor $K_{\tilde \ky_0}+\tilde \kb_0$ is not positive or equivalently we can look at the normalisation $\tilde \ky_0^\nu \to \tilde \ky_0$ and test the curves against the divisor $K_{\tilde \ky_0^\nu}+\nu^*\tilde \kb_0 + D$ where $D$ is the double locus.

Using that $K_{\ky_0}+\kb_0$ was ample we only need to look at the extra curves obtained in the blow-up: the curves on the new component $Q$ and on each of the 4 components of $\ky_0$ one exceptional curve coming from the blow up of a smooth point.

Neglecting the boundary for a moment let $E'$ be any component of the exceptional cycle of lines. Then by computing on any component of $\tilde \ky_0^\nu$ we see that $K_{\tilde \ky_0}$ trivial on $E'$. So we see that $E'$ will not be contrated on the log-canonical model if and only if the boundary divisor $\kb_0$ intersects $E'$ non-trivially. Thus when going to the canonical model, the central components $Q$ will be contracted in one direction to a $\IP^1$ (see Figure 1).
%

The result in the more general case where the surface $\kc$ has $A_n$ singularities can be obtained in a similar manner -- the outcome is the same. 

\begin{figure}[h]\caption{Stable reduction via the normalisation at a degenerate cusp where two components of $\kb_0$ meet.}\label{stabrednormfig}
\centering{\epsfig{file=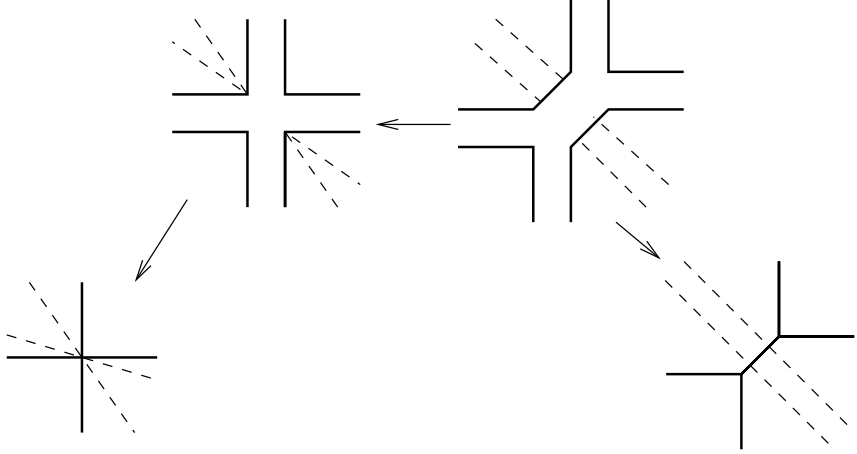}} 
\end{figure}

A posteriori we get the following description of the central fibre:
\begin{prop}\label{centralfibredown}
Let $(\ky_t, \kb_t)_{t\in \Delta^*}$ be a family  as in \ref{assum} Then the stable reduction $(\ky_0, \kb_0)$ of this family can be obtained by the following recipe:
\begin{itemize}
 \item Let $\kc_0$ be the stable degeneration of $\kc_t$ in $\overline{\km}_{g}(H)$, $\nu:\tilde\kc_0\to \kc_0$ the normalisation and $\kp\subset \kc_0$ the  double points.
\item Let $D=\kc_0\times \kp \cup \kp\times \kc_0$ be the double locus of the product $\kc_0\times \kc_0$ and $\kb_0'= \sum_{\phi\in\ks}\lambda_\phi\Gamma_\phi\subset \kc_0\times\kc_0$ be the degeneration of $\kb_t$. 
\item On each component of the normalisation $\tilde \kc_0\times \tilde \kc_0$ blow up the points where  $\bar D=(\nu\times\nu)^*D$ meets the strict  transform of $\kb_0'$.
\item To obtain $\ky_0$ glue everything together along the strict transform of $\bar D$ together with the exceptional curves, $\kb_0$ is then the strict transform of $\kb_0'$.
\end{itemize}
 In short, each degenerate cusp that lies on $\kb_0'$ is replaces by a $\IP^1$ containing two singularities of local type  $\{x_1x_2x_3=0\}\subset \IA^3$.
\end{prop}
The local situation is schematically shown in Figure \ref{stabrednormfig}.

\subsection{The structure of the degeneration}\label{fibs}

The results of the last subsection allow us to describe the degenerations of very simple Galois Kodaira fibrations quite explicitly.
\begin{theo}\label{main}
The moduli space $\MKFbar$ is a union of irreducible components of $\M{st}$.

If $X_0$ is a slc surface corresponding to a point in  $\MKFbar\setminus\MKF$ then
\begin{enumerate}
\item There exists a stable curve $C_0$ and two compact complex surfaces $Y_0$, $Z_0$ fitting in a commutative diagram 
\[
\xymatrix{
X_0\ar[r]^\psi\ar[d]^{\tilde\pi}& Y_0\ar[d]^{\pi}\\
Z_0\ar[r]^{\tilde \psi}& C_0\times C_0
}
\]
such that:
\begin{itemize}
 \item The surface $Y_0$ has slc singularities and $\pi$ is proper and birational. The exceptional locus $E$ of $\pi$ is a disjoint union of $\IP^1$ all of which map to degenerate cusp points in $C_0\times C_0$.The map $\psi$ is a ramified covering. 
\item The surface $Z_0$ does not have slc singularities and the map $\tilde\psi$ is a ramified covering. The map $\tilde\pi$ is proper and birational and contracts the pullback of the exceptional divisor of $\pi$.
\end{itemize}
\item If $n$ denotes the number of irreducible component of $C_0$ then $X_0$ has at least $n^2$ irreducible components. The normalisation of $X_0$ is a union of smooth surfaces, each of which is a (possibly) ramified cover of a product of smooth curves blown up in a finite number of points. 
\item The singular points of $X_0$ are locally isomorphic to one of the following 3 models:
\begin{itemize}
 \item  normal crossing in codimension 1,
\item a degenerate cusp $\spec\IC[x_1, \dots, x_3]/(x_1x_2x_3)$,
\item a degenerate cusp $\spec\IC[x_1, \dots, x_4]/(x_1x_2, x_3x_4)$.
\end{itemize}
 In particular $X_0$ is a local complete intersection and hence Gorenstein.

\end{enumerate}
\end{theo}
\begin{proof}
The surface $X_0$ is in the boundary of $\mathfrak M^{KF}$, so can be obtained as the limit of a family of smooth very simple Galois Kodaira fibrations $\{\psi_t: X_t\to C_t\times C_t\}_{t\in\Delta^*}$ with covering group $G$. We set $C_0$ to be the limit of $C_t$ in the moduli space of stable curves. Then by Lemma \ref{lcmod} $Y_0:=X_0/G$ is the stable reduction of the family of log surfaces $(C_t\times C_t, B_t)$ where $B_t$ is the  branch divisor of $\psi_t$ with appropriate multiplicities; $Y_0$ has been explicitly described in Proposition \ref{centralfibredown} and this description implies the first part of \refenum{i}.

Denoting the complement of the degenerate cusp points in $C_0\times C_0$ with $U$ the $G$-cover $\psi$ induces a $G$-cover over $U$. Then \cite[Lemma 3.2]{ale-pard09} implies that we can complete this to a $G$-cover $\tilde \psi:Z_0\to C_0\times C_0$. By Lemma \ref{quotlem} and Lemma \ref{prodnode} $Z_0$ cannot have slc singularities since the ramification divisor necessarily passes through at least one node. The map $\tilde\pi$ is given by the contraction of $\psi^*E$, the pullback of the exceptional divisor of $\pi$.

Now we prove \refenum{ii}. By Proposition \ref{centralfibredown} components of the  branch divisor on $Y_0$ do not come together in the limit, are smooth when restricted to a component of the normalisation, and do not pass through the degenerate cusps. Therefore also the  ramified coverings of the components of the normalisations are smooth and the number of irreducible components of $X_0$ is at least as big as for $C_0\times C_0$, hence at least $n^2$.

For \refenum{iii} we only have to note that $\psi$ is \'etale near the degenerate cusps of $Y_0$, that have been described in Proposition \ref{centralfibredown}. The ramification meeting the normal crossing locus does not introduce more singularities; again we glue smooth components along a smooth curve.
\end{proof}

\begin{rem}
By definition a very simple Galois Kodaira fibration is in particular a surface fibred over a curve and our description shows that this remains true in the limit. The composition $X_0\to C_0\times C_0 \to C_0$ realises $X_0$ as a fibration over a stable curve with fibres also stable  curves.
The combinatorial type of the fibres, e.g., the number of irreducible components, remains constant on the connected components of $C_0\setminus \mathrm{Sing}(C_0)$ but can change when we pass through a node.
\end{rem}

\section{Deformations and the local structure of the moduli space}\label{def}

The germ of the moduli space of surfaces of general type at a point $[S]$ is a quotient of the versal deformation space $\Def_S$ by $\Aut(S)$ (which may not act faithfully). Unfortunately the latter is very difficult to control so we will only be able to obtain information on $\Def_S$ in some cases.
As a warm-up we will first consider the case of smooth double \'etale Kodaira fibrations before moving on to the stable degenerations.

\subsection{Deformations of $S\in {\gothM}^{KF}$.}\label{defosmooth}
We start in a slightly more general setting. Let $\psi:X\to Y$ be a finite, surjective map  of degree $d$ between smooth surfaces, $B\subset Y$ the branch divisor of $\psi$ and $R\subset X$ the ramification divisor. We assume that both $R$ and $B$ are smooth. Our goal is to obtain information about the deformations of $X$ comparing them to the deformations of the map $\psi$ and these in turn to deformations of $Y$ and the pair $(Y, B)$.

As usual we denote the tangent (resp. obstruction) space to the deformation functor with $\T^1$ (resp. $\T^2$). The infinitesimal automorphisms are denoted by $\T^0$. For example, the infinitesimal deformations of $X$ are classified by $\T_X^1=\Ext^1_X(\Omega_X, \ko_X)=H^1(X, \kt_X)$ and the obstructions lie in $\T_X^2=H^2(X, \kt_X)$. The main tool will be an exact sequence (see \cite{ran89}) where all relevant deformation spaces occur:
\begin{multline}\label{ran}
0\to \T_\psi^0\to \T^0_X\oplus \T^0_Y\to \Ext^0_\psi(\Omega^1_Y, \ko_X)\to  \T_\psi^1\to \T^1_X\oplus \T^1_Y\\
\to \Ext^1_\psi(\Omega^1_Y, \ko_X)\to   \T_\psi^2\to \T^2_X\oplus \T^2_Y\to \Ext^2_\psi(\Omega^1_Y, \ko_X)\to 0.
\end{multline}
The groups $\Ext^k_\psi(\Omega^1_Y, \ko_X)$ can be computed as the limit of two  spectral sequences with $E_2$-term $\Ext^p(L^q\psi^*\Omega_Y^1,\ko_X)$ or $\Ext^p(\Omega_Y^1, R^q\psi_*\ko_X)$ and since the map $\psi$ is flat and finite we have equalities
\begin{align*}
\Ext_\psi^k(\Omega_Y^1, \ko_X)&=\Ext_{\ko_Y}^k(\Omega_Y^1, \psi_*\ko_X)= H^k(Y, \kt_Y\tensor \psi_*\ko_X)\\
&=\Ext_{\ko_X}^k(\psi^*\Omega_Y^1, \ko_X).
\end{align*}

The first characterisation can be used to compare deformations of $X$ and deformations of $\psi$: the push-forward of the structure sheaf of $X$ is locally free of rank $d$ on $Y$ and using the trace map  we can split it as a direct sum 
\beq{Q}\psi_* \ko_X=\ko_Y\oplus \kq,\eeq
the trivial summand corresponding to functions on $X$ which are pullback of functions on $Y$. 
\begin{lem}\label{X=f}
If 
\beq{star} \Ext_Y^i(\Omega_Y, \kq)=0 \text{ for }i=0,1\eeq
then $\T^0_\psi=\T^0_X$,  $\T^1_\psi=\T^1_X$ and $\T^2_\psi\into \T^2_X$, thus $\Def_\psi\to  \Def_X$ is \'etale and every deformation of $X$ is induced by a deformation of the map $\psi$.
\end{lem}
\begin{proof}
The decomposition \refb{Q} induces a decomposition
\[
\Ext_\psi^j(\Omega_Y^1, \ko_X)= H^j(Y, \kt_Y)\oplus\Ext_Y^j(\Omega_Y, \kq) = \T_Y^j\oplus H^j(Y, \kq\tensor\kt_Y)
\]
that combined with \refb{ran} yields an exact sequence
\begin{multline}\label{fX}
 0\to \T^0_\psi\to \T_X^0\to H^0(Y, \kt_Y\tensor \kq)\to \T^1_\psi\to \T_X^1\to H^1(Y, \kt_Y\tensor \kq)\\
\to \T^2_\psi\to \T_X^2\to H^2(Y, \kt_Y\tensor \kq)\to 0.
\end{multline}
By assumption the fourth and the seventh term in the sequence vanish yielding immediately the claimed relationship between $\T^i_\psi$ and $\T^i_X$. The statement on deformations then follows from \cite[Lemma 6.1]{fan-man98}.
\end{proof}

In order to compare further with deformations of $Y$ we have to use the second characterisation of $\Ext_\psi^j(\Omega_Y^1, \ko_X)$. The sheaf of relative differentials $\Omega_{X/Y}$, defined via the exact sequence
\begin{equation}\label{reldiff}
0\to \psi^*\Omega^1_Y\to \Omega^1_X\to  \Omega_{X/Y}\to 0,
\end{equation}
is supported on the ramification divisor $R$. We want to apply the functor $\Hom_X(- , \ko_X)$ to this sequence.
\begin{lem}\label{V^i}
 Let $R=R_1+\dots +R_r$ be a decomposition of $R$ into irreducible components and let $d_i\geq 2$ be the ramification order along $R_i$. Then 
\[\Ext^j_X(\Omega_{X/Y}, \ko_X)\isom \bigoplus_{i=1}^r \bigoplus_{k=2}^{d_i} H^{j-1}(R_i, \ko_{R_i}(kR_i)).\]
In addition, there is a natural inclusion 
\[H^{j-1}(B, \ko_B(B))\into \Ext^j_X(\Omega_{X/Y}, \ko_X)\]
 which in the case where $\psi$ is induced by the action of a group $G$ identifies $H^{j-1}(B, \ko_B(B))$ with $G$-invariant elements in $H^{j-1}(R, \ko_R(\psi^*B))$.
\end{lem}
\begin{proof} The first claim is local along each component of $R$ (recall that we assumed $R$ to be smooth) so that we may assume that $R$ is irreducible and that $\psi:X\to Y$ is cyclic of degree $d$,  ramified along $R$. Let $p$ be a point in $R$ and let $(x_1, x_2)$ be local coordinates centred at $p$ such that locally $R=\{x_1=0\}$ and $\psi(x_1, x_2)=(x^d_1, x_2)$. Thus in a neighbourhood $U$ around $p$ the inclusion $\psi^*\Omega^1_Y\into \Omega^1_X$ in \refb{reldiff} becomes
\[\ko_U d(x_1^d)\oplus \ko_U dx_2\into \ko_U dx_1\oplus \ko_U dx_2.\]
In the second variable, i.e., in the cotangent direction along $R$,  this is an isomorphism. In the first variable we get the structure sheaf of $\ko_{(d-1)R}$ (since $d(x_1^d)=x_1^{d-1}dx_1$) tensored with the conormal bundle $\kn_{R/X}^*=\ko_R(-R)$. As an $\ko_R$-module $\ko_{(d-1)R}$ is generated by powers $1, x_1, x_1^2, \dots, x_1^{d-2}$ of the local equation of $R$, which are also local generators of $ \ko_R(-kR)$. Thus 
\[\Omega_{X/Y}\isom \bigoplus_{k=1}^{d-1} \ko_R(-kR).\] 

To compute the $\Ext$-groups we will use Grothendieck duality for the closed embedding $R\into X$ and thus need the relative dualising sheaf \[\omega_{R/X}= \ko_R(K_R)\tensor {\ko_X(-K_X)\restr R}= \left(\ko_X(K_X+R)\tensor\ko_X(-K_X)\right)\restr R=\ko_R(R).\]
 Then 
\begin{align*}
\Ext^j_X(\Omega_{X/Y}, \ko_X)&=  \Ext^{j-1}_R(\Omega_{X/Y}, \omega_{R/X})\\
&=\Ext^{j-1}_R(\bigoplus_{k=1}^{d-1} \ko_R(-kR), \ko_R(R))\\
&=H^{j-1}(R, \bigoplus_{k=2}^{d} \ko_R(kR)).
\end{align*}

For the second assertion note that $\psi^*B=dR$ and thus $\psi^*\ko_B(B)=\ko_R(dR)$. Taking cohomology and using the above result we get
\begin{multline*}
\Ext^{j+1}_X(\Omega_{X/Y}, \ko_X)\supset H^j(R, \ko_R(dR))=H^j(B, \psi_*\psi^*\ko_B(B))\\=H^j(B, \ko_B(B)\tensor \psi_*\ko_R)\supset H^j(B, \ko_B(B))
\end{multline*}
where in the last step we used the trace map to split off a trivial summand from $\psi_*\ko_R$; this corresponds to splitting off the part that is invariant under the monodromy action which proves the last assertion.
\end{proof}

\begin{prop}\label{abcd} Let $\mathbf V^j:=\Ext^{j+1}_X(\Omega_{X/Y}, \ko_X)/ H^j(B, \ko_B(B))$.
There is a natural exact sequence 
\begin{multline}\label{XYB}
 0\to \T^0_\psi\to \T_{(Y, B)}^0\to \mathbf V^0\to  \T^1_\psi\to \T_{(Y, B)}^1\to 
\mathbf V^1\to  \T^2_f\to \T_{(Y, B)}^2\to 0
\end{multline}
which compares deformations of $\psi$ and deformations of the pair $(Y, B)$.
\end{prop}
\begin{proof}
Both the sequence  \refb{ran} and the long exact sequence obtained by applying $\Hom(-, \ko_X)$ to  \refb{reldiff} feature the spaces $\Ext_\psi^i$ and $\T^i_X$. Comparing them, a simple diagram chase allows  to construct a third long exact sequence
\begin{multline}\label{fY}
0\to  \T^0_\psi\to \T_Y^0\to \Ext^{1}_X(\Omega_{X/Y}, \ko_X)\to \T^1_\psi\to \T_Y^1\\
\to \Ext^{2}_X(\Omega_{X/Y}, \ko_X)
\to \T^2_\psi\to \T_Y^2\to 0.
\end{multline}
Using instead of a direct sum decomposition the short exact sequence
\[0\to H^j(B, \ko_B(B))\to\Ext^{j+1}_X(\Omega_{X/Y}, \ko_X)\to \mathbf V^j\to 0\] 
we can relate \refb{fY} to the deformations sequence of the pair $(Y, B)$  
\begin{multline}\label{YB}
0\to \T_{(Y,B)}^0\to \T^0_Y\to H^0(Y, \ko_B(B)) \to \T_{(Y,B)}^1\to \T^1_Y\\
\to H^1(Y, \ko_B(B))\to \T_{(Y,B)}^2\to \T^2_Y\to 0.
\end{multline}
obtaining the claimed sequence \refb{XYB}.
\end{proof}

Let us now apply the above to very simple Kodaira fibrations.
\begin{theo}\label{defo1}
Let $\psi:X\to Y:=C\times C$ be a very simple Kodaira fibration, $\ks \subset \Aut(C)$ such that the branch divisor $B=\sum_{\sigma\in \ks} \Gamma_\sigma$. Let $H$ be the subgroup of $\Aut(C)$ generated by $\ks$. Assume further condition \refb{star} holds.

Then $\Def_X\isom \Def_{(C, H)}$ is smooth of dimension $\dim H^1(C, \kt_C)^H$ and the germ of the moduli space at $X$  $(\MKF_X,[X])$ is reduced with only finite quotient singularities.

In particular this holds if $\psi$ is a cyclic covering.
\end{theo}
\begin{rem}
For sake of simplicity we stated the Theorem only for very simple Kodaira fibrations but with the necessary change of notation it holds for standard Kodaira fibrations as well. It would be nice if we could improve Theorem \ref{moduliset} so far as to say that the $\mathfrak M_X$ is isomorphic to $\km_{g(C)}(H)$ as a scheme but we do not have sufficient control over the automorphism group of a (standard) Kodaira fibration.

Note that Kodaira fibrations behave quite differently from other ramified cover constructions, say ramified covers of $\IP^2$ or $\IP^1\times \IP^1$. In the latter cases all deformations are induced from deformations of the branch locus in a fixed ambient space, while in our case the branch divisor is rigid in $Y$ and all deformations of $X$ necessarily induce deformations of $Y$.
\end{rem}

\begin{proof}
We have nearly everything in place. The assumptions of Lemma \ref{X=f} are satisfied and thus we can replace $\T^i_\psi$ by $\T^i_X$ in sequence \refb{XYB}. Since $B$ is composed of the graphs of the automorphisms in $\ks$ and these generate the group $H$ we can also replace $\T^i_{(Y, B)}$ by $\T^i_{(C, H)}$. 

For Kodaira fibrations $\mathbf V^0=\Ext^{1}_X(\Omega_{X/Y}, \ko_X)/ H^0(B, \ko_B(B))=0$: let $R_i$ be a component of the ramification divisor mapping to $B_i$. Calculating intersection numbers  we see that $R_i^2 = d B_i^2 = d (2-2g(B_i))<0$ and thus for all $k>0$  the degree of $\ko_R(kR)$ is negative on every component of $R$. In particular, none of these line bundles has global sections and by Lemma \ref{V^i} $\Ext^{1}_X(\Omega_{X/Y}, \ko_X)=0$, thus $\mathbf V^0=0$.

The remaining part of the sequence \refb{XYB} is
\[0\to  \T^1_X\to \T_{(C,H)}^1\to \mathbf V^1\to  \T^2_X\to \T_{(C,H)}^2\to 0\]
and we see that $\T^1_X$ injects into $\T^1_{(C, H)}$.

By \cite{tuff93} we know that $\Def_{(C, H)}$ is smooth of dimension $\dim H^1(C, \kt_C)^H$  and the discussion in Section 6 of \cite{cat-roll06} shows that every deformation of $(C, H)$ (or equivalently $(Y, B)$) induces a deformation of $X$. Thus we have $\T^1_X\isom \T_{(C,H)}^1$, all infinitesimal deformations are unobstructed and $\Def_X$ is smooth of dimension $\dim H^1(C, \kt_C)^H$.

Locally the moduli space $\mathfrak M_X$ is obtained as a quotient of $\Def_X$ by a subgroup of $\Aut(X)$. The latter is finite since $X$ is of general type and thus $\mathfrak M_X$ has only finite quotient singularities.

That cyclic coverings satisfy \refb{star} is the content of the next Lemma.
\end{proof}

\begin{lem}
 If $\psi:X\to Y$ as in Theorem \ref{defo1} is a cyclic covering, then \refb{star} holds.
\end{lem}
\begin{proof} The tangent bundle of $Y=C\times C$ is a direct sum of line bundles $\kt_Y=\mathrm{pr}_1^*\kt_C\oplus\mathrm{pr}_2^* \kt_C$.
If $\psi$ is cyclic of degree $d$ then $\psi_*\ko_X\isom \bigoplus_{i=0}^{d-1}L^{-i}$ for some line bundle $L$ such that $L^d\equiv \ko_Y(B)$, so $\kq=\bigoplus_{j=1}^2\bigoplus_{i=1}^{d-1}L^{-i}\tensor \mathrm{pr}_j^*\kt_C$ and we can treat each of these line bundles separately. A simple calculation of intersection numbers using the Nakai-Moishezon criterion shows that all line bundles $L^i\tensor\mathrm{pr}_j^*\Omega_C$ are ample for $1\leq i \leq d-1$ and by Kodaira vanishing their inverses do not  have cohomology in degree 0 and 1. Thus the same holds for their direct sum and \refb{star} holds in this case.
\end{proof}

We believe that condition \refb{star} holds in many more cases but have no convenient vanishing result to prove this. For example, $\kq$ has been proved to be negative when restricted to curves in $Y$ (see \cite{kebekus-peternell08}) but even if we could show that $\kq^*\tensor \mathrm{pr}^*_j\Omega_C$ is ample this would not yet imply our condition \refb{star}.

\subsection{Deformations of the degenerations}
We can study the deformations of a stable degeneration $X$ in the boundary $\overline{\mathfrak M}^{KF}\setminus{\mathfrak M}^{KF}$ much in the same way as in the smooth case above but several additional difficulties arise. We will need the explict description from Theorem \ref{main} and some additional notation:  there is a diagram
\[\xymatrix{ &X\ar[d]^\psi \\
E \ar@{^{(}->}[r]^r& Y \ar[d]^\pi & B\ar@{_{(}->}[l]_j\ar[d]\\ & C\times C &  B_0\ar@{_{(}->}[l]_{j_0}}\]
where 
\begin{itemize}
 \item $C$ is a stable curve,
\item $B_0=\bigcup_i\Gamma_{\phi_i}$ is a (not necessarily disjoint) union of graphs of automorphisms of $C$,
\item $B$ is the proper transform of $B_0$ and each connected component of $B$ is isomorphic to $C$,
\item $\pi$ is a birational map with exceptional divisor $E=\bigcup_i E_i$ and $E_i\isom \IP^1$,
\item $\psi$ is a ramified covering with Galois group $G$.
\end{itemize}
We denote by $\nu:Y^\nu\to Y$ the normalisation, by $D$ the double locus of $Y$ and by $\bar D$ its pullback by $\nu$. Each component of $\bar D$ is a normal crossing divisor and its singularities map to the degenerate cusps of $Y$.

Note that, since everything is a local complete intersection we still have $T^i_X=\Ext^i(\Omega_X, \ko_X)$ and likewise for $Y$ and $C\times C$.
\subsubsection{Comparing deformations of $X$ and deformations of the pair $(Y, B)$.}
Let $\kq$ and $\mathbf V^i$ be defined as in Section \ref{defosmooth}.
\begin{prop}\label{part1}
 If $\Ext^i(\Omega_Y, \kq)=0$ for $i=0,1$ then we have an exact sequence
\[0\to  \T^1_X\to \T_{(Y,B)}^1\to \mathbf V^1\to  \T^2_X\to \T_{(Y,B)}^2\to 0.\]
The required vanishing holds if the cover $\psi$ is cyclic.
\end{prop}
This can be proved along the same lines as Theorem \ref{defo1} checking that all steps go through in the necessary generality. Both $B$ and the ramification divisor are Cartier divisors, and where they meet the singular locus they locally look like $V(z) \subset \spec \IC[x,y,z]/(x,y)$. A local computation shows that the conclusion of Lemma \ref{V^i} still holds.

Since $B$ is a Cartier divisor in $Y$ and $\Omega_Y$ has no torsion along $B$ we also have the standard exact sequence \refb{YB} for deformations of pairs \cite[Proposition 3.1]{hassett99} and also Proposition \ref{abcd} generalises to this setting.

It remains to check the vanishings $\mathbf V^0=0$, and $\Ext^i(\Omega_Y, \kq)=0$ for $i=0,1$ for cyclic coverings. This can be done componentwise on the normalisation.

\subsubsection{Comparing deformations of the pairs $(Y,B)$ and $(C\times C, B_0)$.}\label{new}
We start with some local computations:
\begin{lem}\label{local}
With $\pi:Y\to C\times C$ as above we have
\begin{enumerate}
 \item $\pi_*\ko_Y=\ko_{C\times C}$ and $R^q\pi_*\ko_Y=0$ for $q>0$.
\item The sequence 
\beq{rel} 0\to \pi^*\Omega_\CC\to \Omega_Y\to \Omega_{Y/\CC}\to 0\eeq
is exact, $\Omega_{Y/\CC}\isom r_*\Omega_E$ and $L^q\pi^*\Omega_\CC=0$ for $q>0$.
\end{enumerate}
\end{lem}
\pf
Both parts can be proved locally: let $A=\IC[x_1, \dots, x_4]/(x_1x_2, x_3x_4)$, $V=\spec A$   and $U=V(zx_1-yx_3, zx_4-yx_2)\subset \IP^1_V$. We get a diagram
\[\xymatrix{ U \ar@{^(->}[r]\ar[d]^\pi & \IA^4\times \IP^1\ar[d]\\ V\ar@{^(->}[r]& \IA^4}\]
The first part of \refenum{i} holds because $\CC$ is Cohen-Macaulay, in particular $S2$. For $q\geq 1$ note that, $V$ being affine, $R^q\pi_*\ko_Y$ is the sheaf associated to the $A$-module $H^q(U, \ko_U)$, which sits in an exact sequence 
\[0=H^q(\IA^4\times\IP^1,\ko_{\IA^4\times\IP^1})\to H^q(U, \ko_U)\to H^{q+1}(\IA^4\times\IP^1, \ki_U).\] Since $\IA^4\times \IP^1$ can be covered by 2 affine patches the \v{C}ech-complex for $\ki_U$ has length 2 and thus $H^{q+1}(\IA^4\times\IP^1, \ki_U)=0$ for $q\geq 1$. So  $H^q(U, \ko_U)=0$ for $q\geq1$ and \refenum{i} follows.

The proof of \refenum{ii} is a straightforward local computation: the sequence
$0\to \ki_V/\ki_V^2\isom \ko_V^{\oplus 2}\overset{\eta}{\to} \Omega_{\IA^4}\restr V \to \Omega_V\to 0$
is a locally free resolution of $\Omega_V$. Taking the pullback with $\pi$ and restricting to the affine subset where (for example) $y\neq 0$ one only needs to check that $\pi^*\eta$ is injective, its cokernel injects in $\Omega_U$ and $\Omega_U\slash\mathrm{coker}(\pi^*\eta)\isom \Omega_{\inverse\pi(0)}$.\qed
\begin{lem}\label{cohom}
 If $L$ is a line bundle on $Y$ then \[\Ext^i_Y(\Omega_{Y/\CC}, L)\isom H^{i-1}(E,   r^*L).\]
\end{lem}
\pf
The inclusion $r:E\into Y$ factors over the normalisation $\nu:Y^\nu\to Y$ by choosing for each component of $E$ one of its preimages in the normalisation, we write  $r=\nu\circ r'$. Then by relative duality both for $\nu$ and  $r'$
\begin{align*} \Ext^i_Y(\Omega_{Y/\CC}, L) &= \Ext^i_Y(r_*\Omega_E, L)= \Ext^i_Y(\nu_*r'_*\Omega_E, L)\\
& = \Ext^i_{Y^\nu}(r'_*\Omega_E, \nu^*L\tensor\omega_\nu) =  \Ext^i_{Y^\nu}(r'_*\Omega_E(\bar D),\nu^*L)\\
& = \Ext^i_{E}(\Omega_E(\bar D), r^*L\tensor \omega_{r'}) =  \Ext^{i-1}_{E}(\Omega_E(\bar D - E), r^*L)\\
& = H^{i-1}(E, \ko_E(E-\bar D -K_E)\tensor  r^*L)= H^{i-1}(E,  r^*L)).
\end{align*}
\qed

\begin{prop}\label{part2}
 All infinitesimal deformations of the pair $(Y,B)$ induce non-trivial deformations of $\CC$ and thus of the pair $(\CC, B_0)$.
\end{prop}
\pf
By Lemma \ref{local} \refenum{i} the spectral sequence  
\[\Ext^p(\Omega_\CC^1, R^q\pi_*\ko_Y)\implies \Ext^{p+q}_\pi(\Omega_{C\times C}, \ko_Y)\]
 degenerates and $\Ext^i_\pi(\Omega_{C\times C}, \ko_Y)=\Ext^i_{C\times C}(\Omega_{C\times C}, \ko_\CC)=\T^i_\CC$. Plugging this in the  sequence \refb{ran} we get isomorphisms $T^i_\pi \isom T^i_Y$ for all $i$.

By Lemma \ref{local} \refenum{ii} the spectral sequence 
\[\Ext^p(L^q\pi^*\Omega_\CC^1,\ko_Y)\implies \Ext^{p+q}_\pi(\Omega_{C\times C}, \ko_Y)\] degenerates as well and we get a long exact sequence as in \refb{fY} for the map $\pi$. Using the isomorphism  $\T^i_\pi\isom \T^i_Y$  there is a diagram with exact rows and columns
\[\xymatrix{ && \Ext^1_Y(\Omega_{Y/\CC}, \ko_Y)\ar[d] \ar[dr]^\delta\\
 0\ar[r] & \T^1_{(Y,B)}\ar[r]\ar[d]& \T^1_\pi\ar[r]\ar[d] & H^1(B, \ko_B(B))\\
& \T^1_{(\CC,B_0)}\ar[r]& \T^1_\CC.}
\]
We claim that the map $\delta$ is injective. Assuming this, every infinitesimal deformation of the pair $(Y,B)$ induces a non-trivial deformation of $\CC$ and, by composition with $\pi$ a deformation of $(\CC, B_0)$.

The map $\delta$ can be factored  $\delta = \beta\circ\alpha$ in the diagram
\[\xymatrix{
\Ext^1_Y(\Omega_{Y/\CC}, \ko_Y)\ar[r]\ar[d]^\alpha & \Ext^1_Y(\Omega_{Y}, \ko_Y)\ar[d]\\
\Ext^1_Y(\Omega_{Y/\CC}, j_*\ko_B)\ar[r]\ar@{=}[d] & \Ext^1_Y(\Omega_{Y}, j_*\ko_B)\ar@{=}[r]&\Ext^1_B(j^*\Omega_{Y}, \ko_B)\ar[d] \\
\Ext^1_B(j^*r_*\Omega_{E}, \ko_B) \ar[rr]^\beta && \Ext^1_B(\ko_B(-B), \ko_B)
} \]

 The kernel of $\alpha$ is a quotient of 
\[\Ext^1_Y(\Omega_{Y/\CC}, \ko_Y(-B))=H^0(E, \ko_E(-B))=0 (\text{Lemma \ref{cohom}})\] so $\alpha$ is injective.

Now we analyse $\beta$: let $D=\sum P_i$ be the Weil-divisor of nodes on $B$, $\nu:\tilde B \to B$ the normalisation and $\bar D=\nu^*D$. Note that $\ko_B(-B)$ is ample on $B$.

Then $j^*r_*\Omega_{E}\isom \ko_D$ and $\beta$ arises by applying $\Ext^i_B(-, \ko_B)$ to the exact sequence
\[ 0\to \ki_D\tensor \ko_B(-B)\to \ko_B(-B)\to \ko_D\to 0.\]
By Serre duality this $\Ext$-sequence is the dual of the long exact sequence in cohomology associated to 
\[ 0\to \ki_D\tensor \omega_B(-B)\to \omega_B(-B)\to \ko_D\to 0\]
and the map $\beta$ is injective if and only if 
\[\beta^\vee: H^0(B, \omega_B(-B))\to H^0(\ko_D)\]
 is surjective if and only if  $H^1(B, \ki_D\tensor \omega_B(-B))=0$ because  $H^1(B, \omega_B(-B))$ vanishes.

Local computations at one node yield 2 exact sequences
\begin{gather*}
 0\to \ki_D\tensor \omega_B(-B)\to \nu_*\nu^* (\ki_D\tensor \omega_B(-B))\to \ko_D^{\oplus2}\to 0,\\
0\to \ko_{\bar D} \to \nu^*(\ki_D\tensor \omega_B(-B))\to \ki_{\bar D}\tensor \nu^*\omega_B(-B)\to 0.
\end{gather*}
Since $H^1(\tilde B, \ko_{\bar D})=0=H^1(\tilde B, \nu^*\omega_B(-B))=H^0(\tilde B, \ko_{\tilde B}(-\bar D)\tensor\nu^*\ko_B(B))^\vee$ we also have $H^1(\tilde B, \nu^*(\ki_D\tensor \omega_B(-B)))=0$.

Combining both sequences to a diagram with exact row and column
{\small
\[\xymatrix{ 0 \ar[d] \\  H^0(\tilde B, \ko_{\bar D}) \ar[d]\ar[dr]^\eta\\
H^0(\tilde B, \nu^* (\ki_D\tensor \omega_B(-B)))\ar[r]^-{\gamma}& H^0(B,\ko_D^{\oplus2}) \ar[r]&
H^1(B, \ki_D\tensor \omega_B(-B)) \ar[r]& 0
}\]}

we see that the composition $\eta$ is an isomorphism, thus $\gamma$ is surjective and $H^1(B, \ki_D\tensor \omega_B(-B))=0$. Therefore $\beta$ and hence $\delta$ are injective which concludes the proof.\qed

We can finally prove the main result of this section where we use the same notation as above.
\begin{theo}\label{main2}
Let ${\gothN}$ be a connected component of $\MKF$  and  $\overline{\gothN}$ its closure in the moduli space of stable surfaces. If for all degenerations $X$ in $\overline{\gothN}\setminus \gothN$ we have $\Ext^i(\Omega_Y, \kq)=0$ for $i=0,1$ then $\overline{\gothN}$ is a connected and irreducible component of the moduli space of stable surfaces.
\end{theo}
\pf Combining Proposition \ref{part1} and Proposition \ref{part2} we see that for every degeneration $X$ we get  $\T^1_X\into \T^1_{\CC, B_0}=\T^1_{(C,H)}$. By our previous study of degenerations this implies $\T^1_X=\T^1_{(C,H)}$ and all deformations of $X$ are in fact degenerations of smooth Kodaira fibrations as described in Theorem \ref{main}. \qed

Since cyclic coverings satisfy the assumptions of the theorem we get
\begin{cor}\label{cyclic}
If $\psi:X\to C\times C$ is a smooth very simple Kodaira fibration such that $\psi$ is  a cyclic covering then $\MKFbar_X$ is a irreducible and connected component of the moduli space of stable surfaces.
\end{cor}

\begin{rem}
 The vanishing conditions we require are not strictly necessary for our result. But if we relax them it might well happen that either there are infinitesimal deformations of $X$ that do not give rise to infinitesimal deformations of $\psi$ or that some infinitesimal deformation of both $X$ and $\psi$ is obstructed for $\psi$ but not for $X$. Does such behaviour actually occur for double Kodaira fibrations or their degenerations?
\end{rem}

%
\end{document}